Educational Production and Optimal Class Size

Thomas H. Foregger
Lucent Technologies, 67 Whippany Road, Whippany, NJ 07981
email: tforegger@comcast.net

ABSTRACT

Lazear (2001) [1] provided a model of a private school with a particular profit function. Using an alternative, related profit function I show that an optimal solution has nearly equal class sizes. I also offer a conjecture about the roots of a certain family of polynomials, which if true, allows one to conclude that if students become less disruptive or teacher cost increases, then the optimal solution for a profit maximizing school is generally, but not always, to use fewer classes. I also show that if the school has $s > 1$ types of students, then the optimal solution will have at most $s - 1$ mixed classes and its bipartite graph is a forest.

**1. Introduction**

A common practice in schools with sufficiently many students in a single grade is to create classes that are are about the same size. Probably this is done for reasons of fairness, so that each teacher has about the same teaching load. One can therefore ask the question: Is this practice optimal in the sense that a profit maximizing private school would do the same thing? It turns out the answer is yes, but a proof of that fact is not easy since it requires solving a nonlinear integer programming problem.

In *Educational Production* Lazear (2001) [1] examines a disruption model of a private school with a particular profit function and states some interesting theorems about the optimal number of classes. In his model there is a probability $p$ that a student is not disrupting the class, a cost $W$ for a teacher, and a number of classes $m$. I use an alternative profit function and prove that if the school is profitable, then the common practice of equalizing class sizes is optimal. The technique relies on Descartes's rule of signs and the arithmetic-geometric mean inequality.

I also consider an extension of the model to allow for $s > 1$ different types of students. In the extended model a student of type $i$ has probability $p_i$ of not disrupting the class. Let $n_{i,j}$ be the number of students of type $i$ in class $j$ ($i = 1, \cdots, s; j = 1, \ldots, m$) for an optimal solution and let $A = (n_{i,j})$. Then the bipartite graph associated with $A$, which has row nodes corresponding to the types of students and column nodes corresponding to the classes is acyclic. Moreover, there are at most $s - 1$ "mixed" classes, i.e. the ones with students of more than one type.

Lazear's model of educational production seems to me to be an important one for the following reason. He creates a profit maximization model for a private school. This allows him to investigate the characteristics of an optimal solution and do static analysis of changes in the model parameters. The conclusions and predicted outcomes seem to explain some of the empirical results in the economics literature about class size and educational outcomes. While the assumptions may seem a bit oversimplified they do allow one to work within a theoretical framework. The test of whether the assumptions are good ones then becomes whether they lead to conclusions and predictions that are supported by empirical studies. Lazear attempted to prove that the class sizes are nearly equal at an optimum, but his proof is flawed. Since the model seems to have some explanatory power, I thought it would be worthwhile to carry out the proof of this result. The fact that classes sizes can be proven to be nearly equal then lends additional credibility to the assumptions of his model.

**2. Specification of the alternative profit function**

Lazear's profit function is

$$f(m, V, W, Z, p) = ZVp^{Z/m} - mW, \qquad (1)$$

where $m > 0$ is the number of classes, $V > 0$ is the value of a unit of learning, $W > 0$ is the cost of a teacher, including rental cost of capital for the classroom, $Z \geq 1$ is the number of students in the school (actually $Z$ is the sum of the class sizes), and $p$ is the probability that a student is not disrupting the class. It is assumed that $p$ is the same for all students and that student behaviors are independent. $Z$ is assumed to be a fixed integer.

Usually a firm has a production function that determines output Q for various inputs. From Q one can then determine revenue, R(Q). In Lazear's model, "educational output" is a "moment of learning in a class of size n" and is equal to $p^n$ for one student and therefore $np^n$ for an entire class. Hence, if class sizes are $n_1, \ldots, n_m$, then total educational output for the school is $\sum_{i=1}^{m} n_i p^{n_i}$. The revenue for this output, assuming the school charges full value, is $V \sum_{i=1}^{m} n_i p^{n_i}$. Since total cost is $mW$, an alternative profit function would be

$$V \sum_{i=1}^{m} n_i p^{n_i} - mW. \qquad (2)$$

Clearly, if $m$ divides $Z$ and $n_i = Z/m$ for all $i$, then Lazear's profit function reduces to (2). Also, if $m = 1$ or $m = Z$ then the two profit functions agree, but they are not the same in general. Since I want to establish that in a profit maximizing school the class sizes are "nearly equal" in the sense that the largest and smallest classes differ by at most 1, and since (1) does not refer to class sizes, I will use the alternative profit function (2). A proof will be given in Theorem 18. The problem of maximizing (2) subject to the constraint that $\sum_{i=1}^{m} n_i = Z$ is a nonlinear integer programming problem with two continuous parameters, $p$ and $W$.

Since this is a long proof I briefly sketch the method here. The basic idea is to let $Q = (k, l)$ be an optimal solution with just two classes and to compare profit at $Q$ with profit at other possible class size vectors. This eliminates $W$ from the problem and one thereby obtains certain polynomials in $p$ that depend on $k$ and $k$ and that must be positive at $p$ in order for $Q$ to be an optimal solution. Consideration of the roots of these polynomials will lead to a contradiction unless $k$ and $l$ differ by at most 1. One can then consider the general case where the solution has $m$ classes. Letting $n_1$ and $n_2$ be the smallest and largest classes, a subschool of $Z' = n_1 + n_2$ students will then have just two classes and so we can conclude that $n_1$ and $n_2$ are nearly equal. The idea of eliminating one cost parameter from the objective function may be useful in other optimization problems.

### 3. Results using the alternative profit function

In this section I use the alternative profit function given by (2) to prove the main result that all class sizes are *nearly equal* at an optimal solution. The proof of this result essentially begins with Proposition 6 and culminates in Theorem 18.

I also assume without loss of generality that $V = 1$, so that $W$ is relative to $V$.

**Proposition 1.** *If $x_1, \ldots, x_m$ are real numbers whose sum is 0, then $\sum_{i=1}^{m} x_i e^{x_i} \geq 0$ with equality if and only if all $x_i$ are 0.*

**Proof.** (due to Robert Israel) Observe that $xe^x \geq x$ for all real $x$, so $\sum_{i=1}^{m} x_i e^{x_i} \geq \sum_{i=1}^{m} x_i = 0$.

**Proposition 2.** *If $x_1, \ldots, x_m$ are real numbers whose sum is 0, then $\sum_{i=1}^{m} e^{x_i}(1 + ax_i)$ is a non-decreasing function of $a$; if not all the $x_i$ are 0, then it is an increasing function of $a$.*

**Proof.** The derivative with respect to $a$ is $\sum_{i=1}^{m} x_i e^{x_i}$ which is nonnegative by Proposition 1.

**Proposition 3.** *For each integer $m \geq 2$ there is a real number $A(m)$ in (-1,0) such that if $a \leq A(m)$ and $x_1, \ldots, x_m$ are real numbers with $\sum_{i=1}^{m} x_i = 0$, and all $x_i \leq -1/a$, then*

$$\sum_{i=1}^{m} e^{x_i}(1 + ax_i) \le m, \tag{3}$$

*and if $a > A(m)$ the inequality (3) does not hold. Also, $A(m) \le A(t)$ whenever $m \ge t$.*

**Proof.** Let $S(m) = \{a | \sum_{i=1}^{m} e^{x_i}(1 + ax_i) \le m$ whenever $\sum_{i=1}^{m} x_i = 0, \ x_i \le -1/a\}$. If $a = 0$ and $x_1, \ldots, x_m$ are not all 0, with $\sum_{i=1}^{m} x_i = 0$, then $\sum_{i=1}^{m} e^{x_i}(1 + ax_i) = \sum_{i=1}^{m} e^{x_i} > m(\prod_{i=1}^{m} e^{x_i})^{1/m} = m$, where the inequality is from the arithmetic-geometric mean inequality. By continuity there exists $\delta > 0$ such that if $-\delta < a < 0$ then $\sum_{i=1}^{m} e^{x_i}(1 + ax_i) > m$; for sufficiently small $a_o$ in $(-\delta, 0)$, $x_i < -1/a_o$ for all $i$. Thus $a_o$ is not in $S(m)$ and $S(m)$ is bounded above by $a_o$. One can easily check that $e^x(1 - x) \le 1$ for all real $x$. Hence if $a = -1$ and $x_i \le 1$ for all $i$, then $\sum_{i=1}^{m} e^{x_i}(1 - x_i) \le m$, so -1 is in $S(m)$. $S(m)$ is therefore a closed, nonempty set which is bounded above, so we can define $A(m)$ to be the maximum element of $S(m)$. If $a \le A(m)$ and $\sum_{i=1}^{m} x_i = 0$, then by Proposition 2, $\sum_{i=1}^{m} e^{x_i}(1 + ax_i) \le \sum_{i=1}^{m} e^{x_i}(1 + A(m)x_i) \le m$. By the definition of $A(m)$, if $a > A(m)$ there exist $x_1, \ldots, x_m$ such that $\sum_{i=1}^{m} x_i = 0$, $x_i \le -1/a$, and $\sum_{i=1}^{m} e^{x_i}(1 + ax_i) > m$.

Claim: $S(m) \subseteq S(t)$ for $m > t$. Let $a$ be in $S(m)$ and suppose $a$ is not in $S(t)$. Then there exist $x_1, \ldots, x_t$, with $x_i \le -1/a$, $\sum_{i=1}^{t} x_i = 0$, such that $\sum_{i=1}^{t} e^{x_i}(1 + ax_i) > t$. Let $x_i = 0$ for $i = t+1, \ldots, m$. Then $\sum_{i=1}^{m} x_i = 0$, $x_i \le -1/a$ and $\sum_{i=1}^{m} e^{x_i}(1 + ax_i) > t + (m - t) = m$, a contradiction to $a$ in $S(m)$.

**Proposition 4.** *Suppose $m \ge 2$ is an integer, $\lambda_1, \ldots, \lambda_m$ are nonnegative reals that sum to 1, and $1 \ge q \ge e^{m/A(m)}$, where $A(m)$ is as in Proposition 3. Then $\sum_{i=1}^{m} \lambda_i q^{\lambda_i} \le q^{1/m}$.*

**Proof.** If $q = 1$ the result is clear so assume $q < 1$. Let $\beta_i = \lambda_i - 1/m$, $a = m/\ln(q)$, and $x_i = \beta_i \ln(q)$. Then

$$x_i \le \frac{-1}{a}, \ \sum_{i=1}^{m} x_i = 0, \ \text{and} \ \sum_{i=1}^{m} \lambda_i q^{\lambda_i} = \sum_{i=1}^{m} (\beta_i + 1/m) q^{\beta_i + 1/m} = \frac{q^{1/m}}{m} \sum_{i=1}^{m} e^{x_i}(1 + ax_i) \le q^{1/m}.$$

The next theorem is not needed for the proof of the main result. I consider a related function $h$ defined as follows:

$$h(m, a, \lambda) = \begin{cases} e^{a/m} - \lambda m & \text{if } m > 0, a < 0, \text{ and } \lambda > 0 \\ 0 & \text{if } m = 0, a < 0, \text{ and } \lambda > 0. \end{cases}$$

Let $\lambda_0 = W/ZV$, $a = Z\ln(p)$. Then

$$ZVp^{Z/m} - Wm = ZV(e^{a/m} - \lambda_o m) = ZVh(m, a, \lambda_o),$$

so for fixed $ZV$, $f$ is a multiple of $h$ and is enough to consider $h$ instead of $f$. Let $n(a, \lambda)$ be the integer solution to the profit maximization problem using $h(m, a, \lambda)$. If there is more than one integer solution, $n(a, \lambda)$ is defined to be the smallest one.

**Theorem 5.** *If the school is profitable using (2) and $n(a, \lambda) > -a$, then the profit at the optimum is less than or equal to the corresponding Lazear profit (1). Equality holds if and only if $p = 1$ or $m$ divides $Z$.*

**Proof.** Let $m = n(a, \lambda)$ be the optimal (integer) number of classes. If $m = 1$ then it is clear that equality holds, so assume $m \geq 2$. By Proposition 3, $-1 < A(m) < 0$ so $m > -a > aA(m)$. Therefore $m/A(m) < a = \ln(p^Z)$. Hence

$$p^Z > e^{m/A(m)}. \tag{4}$$

Let $(n_1, \ldots, n_m)$ be an optimal class size vector and put $\lambda_i = n_i/Z$ for $i = 1, \ldots, m$. Put $q = p^Z$. Then each $\lambda_i > 0$ and $\sum_{i=1}^{m} \lambda_i = 1$. By (4) $1 > q > e^{m/A(m)}$, so that by Proposition 4,

$$\sum_{i=1}^{m} n_i p^{n_i} - mW = Z \sum_{i=1}^{m} \lambda_i q^{\lambda_i} - mW \leq Zq^{1/m} - mW = Zp^{Z/m} - mW.$$

**Proposition 6.** *Let $f(x) = e^x(1 + ax) + e^{-x}(1 - ax)$. (a) If $a \leq -1/2$, then $f(x) \leq 2$ with equality if and only if $x = 0$.*

*(b) If $a \leq -1/2$, then $f'(x) < 0$ if $x > 0$.*

*(c) If $a > -1/2$, then $f(x) > 2$ for small $x > 0$.*

*(d) If $a \leq -1/2$ and $S$ is a discrete set of real numbers, then the maximum value of $f(x)$ for $x$ in $S$ is attained at the element of $S$ with smallest absolute value.*

*(e) $A(2) = -1/2$.*

**Proof.** Expanding $f(x)$ in a power series around $x = 0$ we have

$$f(x) = 2 + 2 \sum_{\substack{k=2 \\ k \text{ even}}} \frac{x^k}{k!} (1 + ak).$$

If $k \geq 2$ then $1 + ak \leq 1 + 2a \leq 0$. Thus each term in the infinite summation is $\leq 0$. From this one can see that (a) and (d) hold. Differentiate the series for $f(x)$ and observe that the coefficients are still $\leq 0$, with the coefficient of $x^3$ negative, so that for $x > 0$, $f'(x) < 0$. This proves (b). To prove (c) notice that $f''(0) = (1+2a)2 > 0$ for $a > -1/2$. So for small $x > 0$, $f''(x) > 0$. Therefore $f'(x)$ is increasing for small positive $x$. But $f'(0) = 0$, so $f'(x) > 0$ for small positive $x$. Since $f(0) = 2$, this shows that $f(x) > 2$ for small positive $x$. (e) follows from (a) and (c).

**Proposition 7.** *Let $g(\lambda) = \lambda q^\lambda + (1-\lambda)q^{1-\lambda}$, where $0 \leq \lambda \leq 1$. Suppose that $1 \geq q \geq e^{-4}$. Then*

*(a) $g(\lambda) \leq q^{1/2}$ with equality if and only if $\lambda = 1/2$ or $q = 1$.*

*(b) if $\lambda$ is restricted to a discrete set in (0,1), the maximum of g is attained at the element of the set that is closest to 1/2.*

**Proof.** If $q$ is 1, the result is clear, so assume $q < 1$. Let $\beta = \lambda - 1/2$, $a = 2/ln(q)$, and $x = \beta \ln(q)$. Then $a \leq -1/2$ and applying Proposition 6, the result follows. The details of applying Proposition 6 are similar to those in the proof of Proposition 4.

**Remark.** *It might be expected that if the school is profitable with profit function (2), then each class is profitable at the optimal solution. However, this is not true. Take $Z = 5, p = .77, W = 1.2$. Then the optimal solution has $m = 2$ and optimal class size vector of (2,3). But $2p^2 - W = -.0142$, $3p^3 - W = .169599$, so the smaller class is not profitable, and the school itself is profitable. The reason for this is that $n_i p^{n_i}$ will not be the same for all $n_i$ unless the $n_i$ are all equal. Hence by choosing W carefully one can make one of the $n_i p^{n_i} - W$ negative and still maintain a profitable school. This result depends on the assumption that Z is fixed. Otherwise, the school would choose to operate with no unprofitable classes.*

**Proposition 8.** *If $1 > d \geq e^{-2}$ and $0 \leq \omega \leq 1$, then $d^\omega(1+\omega) + d^{-\omega}(1-\omega) \leq 2$.*

**Proof.** Put $x = \omega \ln(d)$ and $a = 1/ln(d)$. Then $\ln(d) \geq -2$ implies $a \leq -1/2$. Also, $d^\omega = e^x$, and $\omega = x/ln(d) = ax$. Hence $d^w(1+\omega) + d^{-\omega}(1-\omega) = e^x(1+ax) + e^{-x}(1-ax) \leq 2$ by Proposition 6(a).

**Definition 1.** (a) If $N = (n_1, \ldots, n_m)$ is a class size vector, $V(N, p) = \sum_{i=1}^{m} n_i p^{n_i}$.

(b) For $Z \geq 3$, $\beta_1 = \begin{cases} Z/3 & \text{if } Z = 0(\text{mod } 3) \\ (Z+1)/3 & \text{if } Z = 2(\text{mod } 3) \\ (Z+2)/3 & \text{if } Z = 1(\text{mod } 3), \end{cases}$ $\beta_2 = \begin{cases} Z/3 & \text{if } Z = 0(\text{mod } 3) \\ (Z+1)/3 & \text{if } Z = 2(\text{mod } 3) \\ (Z-1)/3 & \text{if } Z = 1(\text{mod } 3), \end{cases}$

$$\beta_3 = \begin{cases} Z/3 & \text{if } Z = 0 \pmod 3 \\ (Z-2)/3 & \text{if } Z = 2 \pmod 3 \\ (Z-1)/3 & \text{if } Z = 1 \pmod 3. \end{cases}$$

(c) $Q_3(Z) = (\beta_1, \beta_2, \beta_3)$.

(d) $Q_1(Z) = (Z)$.

(e) $Q = (k, l)$, where $k$ and $l$ are positive integers, with $k + l = Z$ and $k \le l$.

(f) $f(p) = 2V(Q, p) - V(Q_3(Z), p) - V(Q_1(Z), p)$.

(g) $\delta_i = Z - 3\beta_i$, $i = 1, 2, 3$.

(h) $\omega = 1/Z$.

(i) $c$ is the unique root of $2x = x^4 + .98$ in $(0,1)$.

Notice that $\beta_3 \le \beta_2 \le \beta_1$ and $\beta_i = (Z - \delta_i)/3$.

**Proposition 9.** *If $Z \ge 6$, then* $\sum_{i=1}^{3} \beta_i p^{\beta_i} \ge .9811 Z p^{Z/3}$.

**Proof.** There is a $C > 0$ such that $p = C^{6/Z}$. Let $R = (1/3) \sum_{i=1}^{3} (1 - \delta_i \omega) C^{-2\delta_i \omega}$. By the arithmetic-geometric mean inequality, and using $\sum_{i=1}^{3} \delta_i = 0$,

$$R^3 \ge \prod_{i=1}^{3} (1 - \delta_i \omega) = 1 + e_2(\delta_1, \delta_2, \delta_3) \omega^2 - \omega^3 \prod_{i=1}^{3} \delta_i,$$

where

$$e_2(\delta_1, \delta_2, \delta_3) = \delta_1 \delta_2 + \delta_1 \delta_3 + \delta_2 \delta_3.$$

Using the actual values for $\delta_i$ we have that the lower bound for $R^3$ is 1 if $Z = 0 \pmod 3$, $1 - 3\omega^2 + 2\omega^3$ if $Z = 1 \pmod 3$, and $1 - 3\omega^2 - 2\omega^3$ if $Z = 2 \pmod 3$. If $Z \ge 8$, then $2\omega^3 + 3\omega^2 = Z^{-2}(2/Z + 3) \le (1/64)(3.25) = .050781$ and $1 + 2\omega^3 - 3\omega^2 \ge 1 - (2\omega^3 + 3\omega^2) \ge .9492$. If $Z = 6$, the lower bound is 1. If $Z = 7$, the lower bound is .9446. Therefore $R \ge .9446^{1/3} = .9811$. It follows that $\sum_{i=1}^{3} \beta_i p^{\beta_i} \ge .9811 Z p^{Z/3}$.

**Proposition 10.** *Suppose $Z \ge 6$. Let $p_1 = c^{6/Z}$. Then $f(p_1) < 0$ and $p_1^Z > e^{-4}$.*

**Proof.** Define $\beta$ by $k = (Z - \beta)/2$, so that $l = (Z + \beta)/2$. Then $\beta \omega < 1$, and $c$ is about $.52922 > .513418 > e^{-2/3}$, so that $p_1^Z > e^{-4}$.

$$f(p_1)c^{-2}/Z = (1-\beta\omega)c^{1-3\beta\omega} + (1+\beta\omega)c^{1+3\beta\omega} - c^4 - (1/3)\sum_{i=1}^{3}(1-\delta_i\omega)c^{-2\delta_i\omega} = L - c^4 - R,$$ where

$$L = (1-\beta\omega)c^{1-3\beta\omega} + (1+\beta\omega)c^{1+3\beta\omega} \text{ and } R = (1/3)\sum_{i=1}^{3}(1-\delta_i\omega)c^{-2\delta_i\omega}.$$ Let $d = c^3$ and $\lambda = d^{\beta\omega}$. Then

$$L = (1-\beta\omega)\frac{c}{\lambda} + (1+\beta\omega)c\lambda = 2c + c(\lambda + \frac{1}{\lambda} - 2 + \beta\omega(\lambda - \frac{1}{\lambda})).$$

Claim: $\lambda + \frac{1}{\lambda} - 2 + \beta\omega(\lambda - \frac{1}{\lambda}) \leq 0$. This is equivalent to $\lambda(1+\beta\omega) + \frac{1}{\lambda}(1-\beta\omega) \leq 2$ or

$$d^{\beta\omega}(1+\beta\omega) + d^{-\omega\beta}(1-\beta\omega) \leq 2, \tag{5}$$

and since $c > e^{-2/3}$ we have $d > e^{-2}$, so (5) is true by Proposition 8. Thus, $L$ is $2c + \gamma$, where $\gamma \leq 0$.

By the proof of Proposition 9 we can write $R = .98 + \gamma_o$, where $\gamma_o > 0$. Hence,

$$f(p_1)c^{-2}/Z = L - c^4 - R = 2c + \gamma - c^4 - .98 - \gamma_o = \gamma - \gamma_o < 0.$$

The fact that $f(p_1) < 0$ and that $p_1^Z > e^{-4}$ turns out to be a crucial result for the proof of the main theorem.

**Proposition 11.** *Suppose $Z \geq 6$. If $(k, Z(mod\ 3)) \neq (\beta_3, 2)$ and $k \geq \beta_3$, then*

*(a) $f$ has exactly one root, $p_o$, in $(0,1)$,*

*(b) 1 is a simple root of $f$,*

*(c) $f(p) > 0$ for $p_o < p < 1$,*

*(d) $f(p) < 0$ for $p_o > p > 0$,*

*(e) $p_o^Z > e^{-4}$.*

**Proof.** The constraints on $k$ and $l$ imply $0 \leq l - k \leq (Z+2)/3$. Hence,

$$f'(1) = 2k^2 + 2l^2 - Z^2 - \sum_{i=1}^{3}\beta_i^2 = (l-k)^2 - \sum_{i=1}^{3}\beta_i^2 \leq ((Z+2)/3)^2 - Z^2/3 = -2/3(Z^2 - 2Z - 2) < 0.$$ By

continuity of $f'$, $f'(p) < 0$ for $p < 1$, $p$ near 1, so $f(p)$ is decreasing for $p$ near 1. But $f(1) = 0$, so $f(p) > 0$ for $p$ near 1. Since $f'(1) \neq 0$, 1 is a simple root of $f$. If $k > \beta_1$, then writing $f$ in order of increasing powers of $p$ we have $f(p) = -\beta_3 p^{\beta_3} - \beta_2 p^{\beta_2} - \beta_1 p^{\beta_1} + 2kp^k + 2lp^l - Zp^Z$. This shows that $f$ has 2 sign changes in its coefficients and for small positive $p$, $f(p) < 0$. Hence, $f$ has at least one root in $(0,1)$. By Descartes's rule of signs, Upsensky (1948) [2], $f$ has 0 or 2 positive roots. But we know that 1 is a simple root of $f$ and there is a root in $(0,1)$. Hence, $f$ has exactly one root in $(0,1)$. Assume that $\beta_3 \leq k \leq \beta_1$. If $Z = 0(mod\ 3)$, then $k = \beta_1 = \beta_2 = \beta_3$ and $f(p) = -\beta_1 p^{\beta_1} + 2lp^l - Zp^Z$ has 2 sign changes in its coefficients, so as before $f$ has exactly one root in (0,1). If $Z = 1 \pmod 3$, and $k = \beta_3 = \beta_2$, then $f(p) = -\beta_1 p^{\beta_1} + 2lp^l - Zp^Z$ has 2 sign changes, so as before $f$ has exactly one root in (0,1). If $Z = 1 \pmod 3$ and $k = \beta_1$ then $f(p) = -2\beta_2 p^{\beta_2} + \beta_1 p^{\beta_1} + 2lp^l - Zp^Z$ has 2 sign changes, so exactly one root in (0,1). If $Z = 2 \pmod 3$, then $k = \beta_1 = \beta_2$ and $f(p) = -\beta_3 p^{\beta_3} + 2lp^l - Zp^Z$ has 2 sign changes and exactly one root in (0,1). It follows that for $0 < p < p_o$, $f(p) < 0$ and for $p_o < p < 1$, $f(p) > 0$. By Proposition 10 we must have $p_o > p_1$. Hence, $p_o^Z > p_1^Z > e^{-4}$.

**Proposition 12.** *Suppose $Z \geq 6$. If $(k, Z(mod\ 3)) = (\beta_3, 2)$, then*

*(a) $f$ has exactly two roots in (0,1), say $p_2$ and $p_3$, and $0 < p_2 < 1/2 < p_3 < 1$,*

*(b) 1 is a simple root of $f$,*

*(c) if $0 < p < p_2$ or $p_3 < p < 1$, then $f(p) > 0$,*

*(d) if $p_2 < p < p_3$, then $f(p) < 0$,*

*(e) $p_3^Z > e^{-4}$.*

**Proof.** We have $f(p) = \beta_3 p^{\beta_3} - 2\beta_1 p^{\beta_1} + 2lp^l - Zp^Z$, so $f$ has 3 sign changes and therefore 1 or 3 positive roots. $f(p) > 0$ for small $p > 0$ since $\beta_3 > 0$. We have $f'(1) = \beta_3^2 - 2\beta_1^2 + 2l^2 - Z^2$ which is $-2(Z+1)(Z-5)/9$, after substituting for $\beta_1$, $\beta_3$, and $l$ and simplifying. Hence 1 is a simple root of $f$, and $f$ is decreasing for $p < 1$, $p$ near 1. Since $f(1) = 0$ it must be that $f(p) > 0$ for $p < 1$, $p$ near 1. We have $f(1/2) = (1/2)^{\beta_3}(-1 + 2l(1/2)^{Z-2\beta_3} - Z(1/2)^{Z-\beta_3}) = (1/2)^{\beta_3}(-1 + (4/3)(Z+1)(1/2)^{(Z+4)/3} - Z(1/2)^{Z-\beta_3})$. One can easily check that $(4/3)(Z+1)(1/2)^{(Z+4)/3}$ is less than or equal to 1 for $Z \geq 5$. Hence $f(1/2) < 0$, and there are roots $p_2$ in (0,1/2) and $p_3$ in (1/2,1). It must be that $p_3 > p_1$, so that $p_3^Z > p_1^Z > e^{-4}$.

Having considered $f$ when $k \geq \beta_3$ we now consider $f$ when $k < \beta_3$.

**Proposition 13.** *Suppose $Z \geq 6$. Let $k_o = Z(3 - 3^{1/2})/6$, if $Z = 0 \pmod 3$ and $k_o = (1/2)Z - (1/6)(3(Z^2+2))^{1/2}$, if $Z \neq 0 \pmod 3$. Assume $k < \beta_3$. Then*

*(a) $k_o < \beta_3$,*

*(b) $f(p) > 0$ for $p > 0$ and $p$ small,*

*(c) if $1 \leq k \leq k_o$, then $f$ has exactly one root in (0,1), say $p_6$, and $p_6 < p_1$, and $f(p) > 0$ if and only if $0 < p < p_6$,*

*(d) if $k_o < k < \beta_3$, then $f$ has exactly two roots in (0,1), say $p_8$, $p_9$ with $1 < p_8 < p_1 < p_9 < 1$.*

**Proof.** (a) It is easy to check that $k_o < \beta_3$.

(b) Since $k < \beta_3$, $f(p) > 0$ for small positive $p$. This completes the proof of (b).

Since $k < \beta_3$, $f$ has 3 sign changes, so 1 or 3 positive roots. By Proposition 10, $f(p_1) < 0$ so $f$ has at least one root in $(0,1)$. Clearly, $f(1) = 0$, so $f$ has 3 positive roots. We calculate

$$f'(1) = 2k^2 + 2l^2 - Z^2 - \sum_{i=1}^{3} \beta_i^2 = (Z - 2k)^2 - \sum_{i=1}^{3} \beta_i^2$$
$$= \begin{cases} (2/3)Z^2 - 4kZ + 4k^2 & \text{if } Z = 0 (\text{mod } 3) \\ (2/3)Z^2 - 4kZ + 4k^2 - (2/3) & \text{if } Z \neq 0 (\text{mod } 3). \end{cases}$$

We view $f'(1)$ as a quadratic polynomial in $k$ and note that $k_o$ is its smallest root. Since $k_o < \beta_3 < Z/2$ and the sum of the roots is $Z$, the other root is greater than $Z/2$. Since $k < \beta_3$, it follows that $f'(1) \neq 0$ if and only if $k \neq k_o$.

(c) If $k < k_0$, then $f'(1) > 0$ so that $f$ is increasing for $p$ near 1. Therefore, $f(p) < 0$ for $p < 1$, $p$ near 1, and $f(p) > 0$ for $p > 1$, $p$ near 1. Since $f(p) < 0$ for large $p$, there is a root greater than 1, so that there is exactly one root in $(0,1)$, and it is less than $p_1$. If $k = k_o$, then $f'(1) = 0$, so that 1 is a double root of $f$, and there is exactly one root in $(0,1)$ and it is less than $p_1$.

(d) If $k > k_o$, then $f'(1) < 0$, so that $f$ is decreasing for $p$ near 1. Therefore, $f(p) > 0$ for $p < 1$, $p$ near 1, and $f(p) < 0$ for $p > 1$, $p$ near 1. It follows that there is a root in $(0, p_1)$ and in $(p_1, 1)$. Since 1 is a root we have accounted for at least 3 positive roots and there are exactly 2 roots in $(0,1)$.

**Proposition 14.** *Suppose* $Z \geq 3$ *and* $k < \beta_3$. *Put* $\tilde{Q}_3 = (k, k, l - k)$ *and* $g(p) = 2V(\tilde{Q}_3, p) - V(\tilde{Q}_3, p) - V(Q_1(Z), p)$. *Then*

*(a) $g$ has exactly one root in $(0,1)$, say $p_4$, and $g(p) > 0$ if and only if $p_4 < p < 1$,*

*(b) $p_1 < p_4$.*

Proposition 14(b) is important for the main result.

**Proof.** (a) $g(p) = -(l - k)p^{l-k} + 2lp^l - Zp^Z$, so that $g$ has 2 sign changes and therefore 0 or 2 positive roots. Since 1 is a root, there are 2 positive roots. For $p$ positive and small, $g(p) < 0$. We have $g'(1) = -(l - k)^2 + 2l^2 - Z^2 = -2k^2 < 0$.

Therefore $g$ is decreasing for $p$ near 1 and $g(p) > 0$ for $p < 1$, $p$ near 1. Hence, there is a root in $(0,1)$, say $p_4$. This must be the only root in $(0,1)$, so that $g(p) > 0$ if and only if $p_4 < p < 1$.

(b) $g(p) = p^{l-k}(-(l - k) + 2lp^k - Zp^{2k}) = p^{l-k}h(p^k)$, where $h(u) = -(l - k) + 2lu - Zu^2$. The quadratic

equation $h(u) = 0$ has one root in $(0,1)$, $(l - k)/Z$, so $p_4 = ((l - k)/Z)^{1/k}$. Let $\alpha = k/Z$, so that $0 < \alpha < 1/3$. To show $p_4 > p_1$, it is sufficient to show that $(1 - 2\alpha)^{1/\alpha} > c^6$. The left side is a decreasing function of $\alpha$ on $(0,1/2)$, and its value at $1/3$ is $1/27 > .0219696 = c^6$.

The following argument shows how to eliminate $W$ from the problem.

Let $Q_m$ be a class size vector for a school of $m$ classes, where $m \geq 3$. If $Q$ is optimal, then

$$V(Q, p) - 2W \geq V(Q_m, p) - mW$$

so that

$$(m - 2)W \geq V(Q_m, p) - V(Q, p). \tag{6}$$

Also,

$$V(Q, p) - 2W > V(Q_1(Z), p) - W$$

so that

$$W < V(Q, p) - V(Q_1(Z), p). \tag{7}$$

Combining (6) and (7) we have

$$V(Q_m, p) - V(Q, p) \leq (m - 2)W < (m - 2)(V(Q, p) - V(Q_1(Z), p))$$

or

$$0 < (m - 1)V(Q, p) - V(Q_m, p) - (m - 2)V(Q_1(Z), p). \tag{8}$$

Notice that when $Q_m = Q_3(Z)$, (8) implies $f(p) > 0$ and when $Q_m = \tilde{Q}_3$, (8) implies $g(p) > 0$.

**Proposition 15.** *Suppose $Z \geq 6$ and $k < \beta_3$. Assume that $Q$ is optimal. Then the components of $Q$ are nearly equal.*

**Proof.** Assume that the components are not nearly equal. By (8) and the definitions of $f$ and $g$, $f(p) > 0$ and $g(p) > 0$. Let $k_o$ be as in Proposition 13. Assume first that $k \leq k_o$. Then $p_4 < p < p_6$ but by Proposition 13(c) and Proposition 14(b), $p_6 < p_1 < p_4$, a contradiction. Next, assume $k_o < k$. Then $p_4 < p < p_8$ (If $p > p_9$, then $p > p_1$ so that $p^Z > e^{-4}$ and the components of $Q$ would be nearly equal.) but by Proposition 13(d) and Proposition 14(b), $p_8 < p_1 < p_4$, a contradiction.

**Proposition 16.** *If $Z \geq 4$, then $(1, Z - 1)$ is not optimal.*

**Proof.** Assume that $(1, Z - 1)$ is optimal. We first consider the case where $p \geq 1/2$. Since $(1, Z - 1)$ is optimal and $Z - 1 \geq 3$, it must be that $p^Z < e^{-4}$, by Proposition 7. If $Z = 4$, then $p < e^{-1} = .3679$, a contradiction to $p \geq 1/2$. If $Z = 5$, then $p < e^{-4/5} = .4493$, a contradiction to $p \geq 1/2$. Hence we can assume $Z \geq 6$. Suppose first that $Z$ is odd, so that $Z \geq 7$. Let $m = (Z+1)/2$ and $Q_m = (1, 2, \ldots, 2)$. Then the right side of (8) is

$$h(p) := (1/2)(Z-1)(p + (Z-1)p^{Z-1}) - p - (Z-1)p^2 - (1/2)(Z-3)Zp^Z$$

$$= (1/2)(Z-3)p - (Z-1)p^2 + (1/2)(Z-1)^2 p^{Z-1} - (1/2)Z(Z-3)p^Z$$

which has 3 sign changes. So there are 1 or 3 positive roots. For $p > 0$ and small, $h(p)$ is positive. One can easily check that $h'(1) = 0$, so that 1 is a double root of $h$. We have $h(1/2) = -1/2 + (1/2)^{Z+1}(Z^2 + 7Z + 2) < 0$ since $Z \geq 7$. Hence the third root is in $(0, 1/2)$ and $h(p) < 0$ for $p$ in $(1/2, 1)$. This is a contradiction to (8). Second, assume that $Z$ is even. Let $m = Z/2$ and $Q_m = (2, , \ldots, 2)$. Then the right side of (8) is

$$h(p) := (1/2)(Z-2)(p + (Z-1)p^{Z-1}) - Zp^2 - (1/2)(Z-4)Zp^Z$$

$$= (1/2)(Z-2)p - Zp^2 + (1/2)(Z-2)(Z-1)p^{Z-1} - (1/2)(Z-4)Zp^Z$$

which has 3 sign changes, so again there are 1 or 3 positive roots. We have $h(p) > 0$ for small positive $p$ and $h'(1) = Z - 1$, so that $h'(1) > 0$ and 1 is a simple root of $h$, and $h$ is increasing for $p < 1$, $p$ near 1. Since $h(1) = 0$ it must be that $h(p) < 0$ for $p < 1$, $p$ near 1. We have $h(1/2) = -1/2 + (1/2)^{Z+1}(Z^2 - 2Z + 4) < 0$, so that there is a root in $(0, 1/2)$. Since $h(p) < 0$ for large $p$ there is a third root in $(1, \infty)$. Hence $h(p) < 0$ for $p$ in $(1/2, 1)$. This is a contradiction to (8). This completes the case where $p \geq 1/2$. Now assume $p < 1/2$. Since $(1, Z-1)$ is optimal, we have $(Z-1)p^{Z-1} + p - 2W > Zp^Z - W$ and $(Z-1)p^{Z-1} + p - 2W \geq (Z-2)p^{Z-2} + 2p - 3W$. It follows from these two inequalities that $(Z-1)p^{Z-1} + p - 2W > \dfrac{Zp^Z + (Z-2)p^{Z-2}}{2} + p - 2W$, so that

$$(Z-1)p^{Z-1} > \dfrac{Zp^Z + (Z-2)p^{Z-2}}{2} \text{ or}$$

$$h(p) := -Zp^2 + 2(Z-1)p - (Z-2) > 0.$$

$h$ has 2 sign changes, so 0 or 2 positive roots. Since 1 is a root, there are 2 positive roots. $h(p) < 0$ for small positive $p$ and $h'(1) = -2$ so $h(p) > 0$ for $p < 1$, $p$ near 1. Hence $h$ has a unique root in $(0, 1)$. $h(1/2) = -Z/4 + 1 \leq 0$ since $Z \geq 4$. Hence $h$ is negative on $(0, 1/2)$, which is a contradiction.

**Proposition 17.** *Suppose* $Z \geq 6$. *If the optimal number of classes is 2, then the class sizes are nearly*

*equal.*

**Proof.** Suppose $Q$ is an optimal class size vector. If $k < \beta_3$, we can apply Proposition 15, so assume $k \geq \beta_3$. From (8) with $Q_3 = Q_3(Z)$ it follows that $f(p) > 0$. If $(k, Z(\text{mod } 3)) \neq (\beta_3, 2)$, then by Proposition 11, $p > p_o$, so that $p^Z > p_o^Z > e^{-4}$. If $(k, Z(\text{mod } 3)) = (\beta_3, 2)$ then by Proposition 12, $p > p_3$, so that $p^Z > p_3^Z > e^{-4}$. It follows from Proposition 7 that $k$ and $l$ are nearly equal.

The following theorem is the main result of this paper.

**Theorem 18.** *All class sizes are nearly equal at an optimal solution.*

**Proof.** Let $m$ be the optimal number of classes and let $(n_1, \ldots, n_m)$ be an optimal class size vector. We can assume that $2 \leq m < Z$. Let $n_1$ be the smallest class, $n_m$ the largest class, and $Z' = n_1 + n_m$. Then $(n_1, n_m)$ is optimal for a school of $Z'$ students, and $Z' \geq 2$. If $Z' = 2$ or $Z' = 3$ then clearly we have $n_1$ and $n_m$ nearly equal. If $Z' = 4$ or $Z' = 5$ we can apply Proposition 16 and if $Z' \geq 6$ we can apply Proposition 17 to conclude $n_1$ and $n_m$ are nearly equal. Since the largest and smallest class sizes differ by at most 1, the class sizes are nearly equal.

**Remark.** *It should be noted that if one assumes that only a predetermined number of classes will run, then it may not be optimal to have them nearly equal. For example, suppose that the school has $Z = 100$ students and for whatever reason, only two classes will be run. Let $p = .95$. Then a class size vector of (23,77) has revenue of 8.55, while a vector of (50,50) has only 7.69. The optimum solution is (23,77). The reason for this is that the function $g(\lambda) = \lambda q^\lambda + (1 - \lambda) q^{1-\lambda}$ has its maximum on (0,1) at a value of $\lambda \neq 1/2$, when $q < e^{-4}$. So if $q = p^Z < e^{-4}$, the optimum class size vector will have unequal components.*

### 4. Changes in p or W

The next major result is Theorem 30 which characterizes what happens to the optimal number of classes as $p$ or $W$ increases. I am unable to give a proof of Theorem 30 which does not depend on Conjecture A, which is stated below. Propositions 21 to 29 are in support of Theorem 30.

**Definition 2.** (a) $M_o = Z/2$ if $Z$ is even and $(Z + 1)/2$ if $Z$ is odd.

(b) For each integer $k$, $1 \leq k \leq M_o + 1$ let $Z = kq_k + r_k$, where $q_k$ and $r_k$ are integers and $0 \leq r_k < k$.

(c) $Q_k = (q_k, q_k, \ldots, q_k, q_k + 1, \ldots, q_k + 1)$, where there are $k$ components and $r_k$ of them are equal to $q_k + 1$.

(d) $f_k(p) = V(Q_k, p) - V(Q_{k-1}, p)$ $(2 \leq k \leq M_o + 1)$.

**Proposition 19.** *Let $q_k$ and $r_k$ be as in Definiton 2.*

(a) *If $k \geq 2$, then $q_k \leq q_{k-1}$.*

(b) *If $k \geq 2$ and $q_k = q_{k-1}$, then $r_{k-1} > r_k$.*

**Proof.** (a) $q_k = [Z/k] \leq [Z/(k-1)] = q_{k-1}$, where [] denotes the floor function.

(b) We have $kq_k + r_k = (k-1)q_{k-1} + r_{k-1}$ so that $0 < q_k = r_{k-1} - r_k$.

**Proposition 20.** *Let $q_k$, $r_k$, and $f_k$ be as in Definition 2. Define $S(k) = (k - r_k)q_k^2 + r_k(q_k + 1)^2$. Then for $k > 1$*

*(a) $S(k-1) + S(k+1) \geq 2S(k)$, with equality if and only if*

   i.  $q_k = q_{k+1} = q_{k-1}$, *or*

   ii. $k = Z/d + 1$, *where $d$ divides $Z$, $d \geq 2$, $Z \geq 2d(d-1)$, $q_{k+1} = q_k = d - 1$ and $q_{k-1} = d$.*

*(b) If $j > k$, then $f'_j(1) \geq f'_k(1)$.*

**Proof.** (a) (This is due to Peter Winkler) One can easily check that $S(k) = Zq_k + r_k(q_k + 1)$. For integer $j$, define

$$S_j(k) = (q_k - j)Z + (r_k + kj)(q_k - j + 1).$$

Then $S_j(k) = Zq_k + r_k(q_k + 1) - jZ + jkq_k + jr_k - 2jr_k - kj(j-1) = S(k) - 2jr_k - kj(j-1)$. If $j > 0$, then $2jr_k$ and $kj(j-1)$ are both non-negative, so that $S_j(k) \leq S(k)$. If $j < 0$, then $kj(j-1) \geq k2|j| > 2|j|r_k$, so that $-2jr_k - kj(j-1) < 0$ and $S_j(k) < S(k)$ when $j$ is negative. Now write $q_{k-1} = q_k + j$ and $q_{k+1} = q_k + i$, where $j \geq 0$, and $i \leq 0$ by Proposition 19(a). Observe that $S(k-1) + S(k+1) \geq S_j(k-1) + S_i(k+1) = 2S(k)$. If equality holds, then $i = j = 0$, or $i = 0, j = 1$ and $r_{k-1} = 0$, since one must have $-2jr_{k-1} - j(j-1)(k-1) = 0$. When $i = j = 0$ we have $q_k = q_{k+1} = q_{k-1}$. When $i = 0$ and $j = 1$ and $r_{k-1} = 0$ we have $q_{k+1} = q_k$ and $k = Z/d + 1$ for some integer $d$ which divides $Z$, and $q_{k-1} = q_k + 1$. We have $r_k = Z - kq_k = r_{k-1} + k - d = k - d$, so that $0 \leq q_{k+1} = r_k - q_k = (k - d) - (d - 1) = Z/d + 2 - 2d$, so $Z \geq 2d(d-1)$. Since $k < Z$ it follows that $d \geq 2$.

(b) From (a) we have $f'_{k+1}(1) = V'(Q_{k+1}, 1) - V'(Q_k, 1) = S(k+1) - S(k) \geq S(k) - S(k-1) = f'_k(1)$. Hence, by induction on $j$, $f'_j(1) \geq f'_k(1)$.

**Proposition 21.** *Let $f_k$ be as in Definition 2 and $S(k)$ be as in Proposition 20. The following are equivalent:*

*(a) The equality conditions in Proposition 20(a).*

*(b) $f_{k+1} = f_k$.*

*(c) $f'_{k+1}(1) = f'_k(1)$.*

*(d) $S(k-1) + S(k+1) = 2S(k)$.*

**Proof.** Suppose (a) holds. If (i) holds, then $q_k = q_{k-1} = q_{k+1} = q$, say. It is easy to check that $r_{k-1} + r_{k+1} = 2r_k$ in this case. Hence, $V(Q_{k-1}, p) + V(Q_{k+1}, p)$
$= (2k - r_{k-1} - r_{k+1})qp^q + (r_{k-1} + r_{k+1})(q+1)p^{q+1} = 2((k-r_k)qp^q + r_k(q+1)p^{q+1}) = 2V(Q_k, p)$, so it follows that $f_{k+1} = f_k$.

If (ii) holds, then it is easy to check that $r_{k-1} = 0$, $r_k = k - d$, and $r_{k+1} = k + 1 - 2d$. From this it follows again that $V(Q_{k-1}, p) + V(Q_{k+1}, p) = 2V(Q_k, p)$, and that $f_{k+1} = f_k$.

That (b) implies (c) is trivial. If (c) holds then

$$0 = f'_{k+1}(1) - f'_k(1) = V'(Q_{k+1}, 1) - V'(Q_k, 1) - (V'(Q_k, 1) - V'(Q_{k-1}, 1))$$
$$= S(k+1) - S(k) - (S(k) - S(k-1))$$
$$= S(k+1) + S(k-1) - 2S(k),$$

so that (d) holds. If (d) holds, then (a) holds by Proposition 20(a).

**Proposition 22.** *Let $f_k$ be as in Definition 2. Then*

*(a) If $f_i \neq f_j$ there is a unique root $p_{i,j}$ of $f_i - f_j$ in $(0,1)$. Moreover, $p_{i,j}$ is a simple root.*

*(b) If $0 < p < 1$, and $j > i$, then $f_j(p) - f_i(p) > 0$ if and only if $0 < p < p_{i,j}$.*

*(c) If $i \geq 2$ there is a unique $p_{i,1}$ in $(0,1)$ such that $f_i(p_{i,1}) = V(Q_1, p_{i,1})$.*

*(d) There is an $s_k$ in $(0,1)$ such that $f_k(p)$ is increasing for $p < s_k$ and decreasing for $p > s_k$.*

**Proof.** (a,b) One can easily check that the smallest degree term of $f_k$ is

$$(k - r_k)q_k p^{q_k} \text{ if } q_k < q_{k-1},$$
$$(1 + q_k)q_k p^{q_k} \text{ if } q_k = q_{k-1}.$$

Now let $F = f_j - f_i$ where $j > i \geq 2$.

Claim: $F(p) > 0$ for small positive $p$.

We have $q_j \leq q_{j-1} \leq q_i \leq q_{i-1}$. We have two cases to consider.

Case A. $q_j < q_i$. In this case $f_j$ has the smallest degree term, and its coefficient is positive, so the smallest degree term of $F$ has a positive coefficient and $F(p) > 0$ for small positive $p$.

Case B. $q_j = q_i$. In this case $q_j = q_{j-1}$. If $q_i = q_{i-1}$, then we have $f_i = f_j$ by Proposition 21, a contradiction. Hence, $q_i < q_{i-1}$. Let $q_{i-1} = q_i + b$, where $b \geq 1$. The smallest degree term of $F$ is then

$(1 + q_i)q_i p^{q_i} - (i - r_i)q_i p^{q_i} = q_i p^{q_i}(1 + r_{i-1} - i + b(i-1)) \geq q_i p^{q_i}(1 + r_{i-1} - i + i - 1) = q_i p^{q_i} r_{i-1} \geq 0$, with

equality if and only if $r_{i-1} = 0$ and $b = 1$. But if equality holds, then $f_i = f_{i+1} = \cdots = f_j$, a contradiction. Hence $F(p) > 0$ for small positive $p$.

Since $j > i$ and $f_j \neq f_i$ we have $F'(1) = f'_j(1) - f'_i(1) > 0$ by Proposition 21 and Proposition 20(b). Therefore, $F(p) < 0$ for $p < 1$, $p$ near 1. The two conditions on $F$ give that $F$ has a root in $(0,1)$. The number of sign changes of $F$ is 2 so there is exactly one root in $(0,1)$. Denote this root by $p_{i,j}$. If $0 < p < 1$, then $F(p) > 0$ if and only if $0 < p < p_{i,j}$.

(c) Let $F(p) = f_i(p) - V(Q_1, p)$. Then the smallest degree term of $F$ is $Aq_i p^{q_i}$, where $A > 0$. Also, $F(1) = f_i(1) - Z = -Z$. Hence $F(p) > 0$ for small positive $p$, and $F(1) < 0$. Therefore, $F$ has a unique root, $p_{i,1}$ in $(0,1)$.

(d) By Proposition 19(a) $q_k \leq q_{k-1}$. We have

$$f_k(p) = V(Q_k, p) - V(Q_{k-1}, p) = Ap^{q_k} + Bp^{q_k+1} - Cp^{q_{k-1}} - Dp^{q_{k-1}+1} \quad (2 \leq k \leq M_o + 1)$$

for positive integers $A$ and $C$ and nonnegative integers $B$ and $D$. By Proposition 19(a) there are two cases: (1) $q_k < q_{k-1}$ and (2) $q_k = q_{k-1}$. In case (1) we have $q_k \leq q_{k-1} - 1$, and

$$f'_k(p) = Aq_k p^{q_k-1} + B(q_k + 1)p^{q_k} - Cq_{k-1}p^{q_{k-1}-1} - D(q_{k-1} + 1)p^{q_{k-1}}. \tag{9}$$

We have

$$q_k - 1 < q_k \leq q_{k-1} - 1 < q_{k-1}$$

so that (9) is written in order of increasing powers of $p$. Therefore $f'_k(p)$ has at most one sign change in its coefficients and at most one positive root. Also, $f'_k(p) > 0$ for small positive $p$ and $f'_k(p) < 0$ for large $p$, so that $f'_k$ has exactly one positive root. In case (2) we have

$$f'_k(p) = Aq_k p^{q_k-1} - Cq_{k-1}p^{q_{k-1}-1} + B(q_k + 1)p^{q_k} - D(q_{k-1} + 1)p^{q_{k-1}}$$

$$= q_k(A - C)p^{q_k-1} + (q_k + 1)(B - D)p^{q_k},$$

so that $f'_k$ at most one sign change and at most one positive root. We use the following:

$$A = q_k(k - r_k),$$
$$B = r_k(q_k + 1),$$
$$C = q_{k-1}(k - 1 - r_{k-1}),$$
$$D = r_{k-1}(q_{k-1} + 1).$$

By Proposition 19(b), $r_{k-1} > r_k$. Hence, $A - C = q_k(1 + r_{k-1} - r_k) \geq 2q_k > 0$ and $B - D = (q_k + 1)(r_k - r_{k-1}) < 0$, and again $f_k'(p) > 0$ for small positive $p$ and $f_k'(p) < 0$ for large $p$, so that $f_k'$ has exactly one positive root. This completes the proof of case (2). Since $f_k(0) = f_k(1) = 0$, this root must be in (0,1). Let $f_k'(s_k) = 0$, so that $f_k(p)$ is increasing for $p < s_k$ and decreasing for $p > s_k$. This completes the proof of the Proposition.

**Definition 3.** $R(m) = \{(p, W) \mid \text{the optimal solution has } m \text{ classes and the school is profitable}\}$.

**Remark.** *It is a convention that the given point $(p, W)$ belongs to the set $R(m)$ with least $m$.*

The next proposition is a technical result that is only used in Proposition 25.

**Proposition 23.** *Suppose that $p \geq 1/2$, $W < 2p(1 - p)$, and $Z = bm$, where $b \geq 2$ and $m$ are integers. Then $Zp^{Z/m} - Wm < (p - W)Z$.*

**Proof.** Define a function $q(p) = (1 - p^{b-1})/(1 - p)$. Since $m < Z$ we have

$$Zp^b - Wm - (pZ - WZ) = Z(p^b - p) + W(Z - m) < -Z(p - p^b) + 2p(1 - p)(Z - m)$$

$$= -Zp(1 - p)q(p) + 2p(1 - p)(Z - m) = mp(1 - p)(b(2 - q(p)) - 2).$$

Since $p \geq 1/2$, we have $2 - q(p) \leq 2 - q(1/2) = (1/2)^{b-2}$. Also, $2^{b-1} \geq b$. Hence $b(2 - q(p)) \leq b(1/2)^{b-2} \leq 2^{b-1}(1/2)^{b-2} = 2$. Therefore, $Zp^b - Wm < pZ - WZ$.

**Proposition 24.** *If $n \geq 1$ is an integer and $0 < a \leq 1/2$, then $na^n \leq a$.*

**Proof.** This is clear if $n = 1$. If $n \geq 2$, then $a \leq 1/2 \leq (n - 1)/n$, so that $na \leq n - 1$ and $na^n \leq (n - 1)a^{n-1}$, and the result follows by induction on $n$.

**Proposition 25.** $R(Z) = \{(p, W) \mid W < 2p(1 - p) \text{ and } W < p\}$.

**Proof.** Let $(n_1, \ldots, n_m)$ be an optimal solution for (2). Assume first that $p \geq 1/2$ and $W < 2p(1 - p)$. Then clearly the school is profitable using all classes of size 1. If $m = Z$, then clearly $(p, W)$ is in $R(Z)$, so assume $m < Z$. For each integer $k \geq 1$ let $S_k = \{i \mid n_i = k\}$ and $m_k$ be the number of elements in $S_k$. Then $m_1 \neq m$ because $m < Z$. Let $Z_k = km_k$. Then $\sum_{k=1}^{\infty} Z_k = Z$ and using Proposition 23 we have

$$\sum_{i=1}^{m} n_i p^{n_i} - mW = \sum_{k=1}^{\infty} \sum_{i \varepsilon S_k} n_i p^{n_i} - mW = \sum_{k=1}^{\infty} \sum_{i \varepsilon S_k} kp^k - mW = \sum_{k=1}^{\infty} (m_k kp^k - Wm_k)$$

$$= \sum_{k=1}^{\infty} (Z_k p^{Z_k/m_k} - Wm_k) < \sum_{k=1}^{\infty} (pZ_k - WZ_k) \text{ (since } m_1 \neq m)$$

$$= (p - W)Z = \text{ profit for } Z \text{ classes of size 1 each.}$$

Hence $(p, W)$ is in $R(Z)$. Now suppose $p < 1/2$ and $W < p$. Then by Proposition 24

$$\sum_{i=1}^{m} (n_i p^{n_i} - W) \leq \sum_{i=1}^{m} (p - W) = m(p - W) \leq Z(p - W) = \text{ profit for } Z \text{ classes of size 1 each,}$$

and $Z(p - W) > 0$ so the school is profitable. If equality holds, then each $n_i = 1$ and $m = Z$. Hence $(p, W)$ is in $R(Z)$. Suppose $(p, W)$ is in $R(Z)$. If $Z$ is even, the profit from $Z$ classes of size 1 exceeds the profit from $Z/2$ classes each of size 2. Thus, $Zp - ZW > (Z/2)(2p^2) - WZ/2 = Zp^2 - WZ/2$ so $2p(1 - p) > W$. If $Z$ is odd, a similar argument with one class of size 1 and $(Z - 1)/2$ classes of size 2 also shows that $2p(1 - p) > W$. Since the school is profitable, $p > W$.

**Proposition 26.** *If $W \geq 2p(1 - p)$, then an optimal solution has at most one singleton class.*

**Proof.** Suppose there are two singleton classes in an optimal solution. Let $m$ be the number of classes and let $(1, 1, n_1, \ldots, n_{m-2})$ be an optimal class size vector. Then $(2, n_1, \ldots, n_{m-2})$ has $m - 1$ classes, so has strictly smaller profit. Therefore, $2p + \sum_{i=1}^{m-2} n_i p^{n_i} - mW > 2p^2 + \sum_{i=1}^{m-2} n_i p^{n_i} - W(m - 1)$, which implies $W < 2p(1 - p)$, a contradiction.

**Proposition 27.** *If $p \leq 1/2$ and the school is profitable, then $p > W$ and the optimal class sizes are 1.*

**Proof.** Let $m$ be the optimal number of classes and let $(n_1, \ldots, n_m)$ be an optimal class size vector. Then using Proposition 24

$$0 < \sum_{i=1}^{m} n_i p^{n_i} - mW = \sum_{i=1}^{m} (n_i p^{n_i} - p) + (p - W)m \leq (p - W)m,$$

so $p > W$. But also, $(p - W)m \leq (p - W)Z = $ profit from a school of $Z$ classes each of size 1.

**Theorem 28.** *For a profitable school of $Z$ students the optimal number of classes $m$ satisfies $m = Z$, or $m \leq Z/2$, $Z$ even or $m \leq (Z + 1)/2$, $Z$ odd.*

**Proof.** Suppose first that $W < 2p(1 - p)$. If $p \leq 1/2$, then by Proposition 27, $p > W$ and all classes have size 1, so $m = Z$. If $p > 1/2$ then clearly $W < p$, so $(p, W)$ belongs to $R(Z)$, by Proposition 25. So we can assume $W \geq 2p(1 - p)$. By Proposition 26, there is at most one singleton class. Let the optimal class size vector be $(n_1, \ldots, n_m)$. If each $n_i \geq 2$ then $Z = \sum_{i=1}^{m} n_i \geq 2m$, so $m \leq Z/2$. If $n_1 = 1$, then

$n_i \geq 2, i = 2, \ldots, m$, so $Z = \sum_{i=1}^{m} n_i \geq 1 + 2(m - 1)$ and $m \leq (Z + 1)/2$. If $Z$ is even, then $m \leq Z/2$.

Numerical evidence supports the following conjecture:

**Conjecture A.** If $i > k + 1$ then $p_{i,k} \leq p_{k,k+1}$. If $i < k$, then $p_{i,k+1} \geq p_{k,k+1}$.

Based on this conjecture I can show the following theorem:

**Proposition 29.** *(a) If $p_{i,k}$ and $p_{i,j}$ are defined, and $k \geq j$, then $p_{i,k} \leq p_{i,j}$.*

*(b) If $f_i \neq f_2$, then $p_{2,i} \leq p_{2,1}$.*

**Proof.** (a) It is enough to show that

(i) $p_{i,k+1} \leq p_{i,k}$ and (ii) $p_{i+1,i} \leq p_{i,i-1}$.

Assuming (i) I can show (ii): $p_{i+1,i} \leq p_{i+1,i-1} \leq p_{i,i-1}$.

If $f_k = f_{k+1}$ there is nothing to show, so assume that $f_k \neq f_{k+1}$. Then $p_{k,k+1}$ exists.

Claim: $p_{i,k} \geq p_{i,k+1}$ whenever $i > k + 1$.

By Conjecture A, $p_{k,k+1} \geq p_{i,k}$, so that $(f_{k+1} - f_k)(p_{i,k}) \geq 0$. Hence,

$(f_i - f_{k+1})(p_{i,k}) \leq (f_i - f_k)(p_{i,k}) = 0$, so that $p_{i,k} \geq p_{i,k+1}$.

Now assume that $i < k + 1$. Since $f_i \neq f_k$, we have $i < k$. By Conjecture A, $p_{k,k+1} \leq p_{i,k+1}$, so that

$(f_{k+1} - f_k)(p_{i,k+1}) \leq 0$. Hence, $0 = (f_{k+1} - f_i)(p_{i,k+1}) \leq (f_k - f_i)(p_{i,k+1})$, so that $p_{i,k+1} \leq p_{i,k}$.

(b) Since $f_i \neq f_2$, we have $i > 2$ and $Z \geq 5$. By direct computation one can check that $p_{2,3} < p_{2,1}$ if $Z = 5$ or

$Z = 6$. So we assume $Z \geq 7$. By (a) $p_{2,i} \leq p_{2,3}$. Notice that $Q_2 = Q_2(Z)$. Put $\omega = 1/Z$ and $p = c^\omega$, where

$c = 12/49$. If $Z$ is odd, then $V(Q_2, p) - 2V(Q_1, p) = Z[(1/2)(1 - \omega)c^{(1-\omega)/2} + (1/2)(1 + \omega)c^{(1+\omega)/2} - 2c]$

$> Z[(1 - \omega^2)c]^{1/2} - 2Zc \geq Zc^{1/2}(4/7)\sqrt{3} - 2Zc$ (since $Z \geq 7$) $= 0$. Hence $f_2(p) > V(Q_1, p)$. If $Z$ is even, then

$V(Q_2, p) - 2V(Q_1, p) = Zp^{Z/2} - 2Zp^Z = Z(c^{1/2} - 2c) = Z((2/7)\sqrt{3} - 24/49) > 0$, so again $f_2(p) > V(Q_1, p)$.

Thus, in both cases, $p_{2,1} \geq c^{1/Z}$. Note that $Q_3 = Q_3(Z)$ and write $Q_2(Z) = ((Z - \beta)/2, (Z + \beta)/2)$, where $\beta$

is 0 or 1. Put $a = c^{1/6}$. Using Definition 2, we have $f_2(p) - f_3(p) = 2V(Q_2, p) - V(Q_1, p) - V(Q_3, p)$

$= (Z - \beta)p^{(Z-\beta)/2} + (Z + \beta)p^{(Z+\beta)/2} - Zp^Z - \sum_{i=1}^{3}(1/3)(Z - \delta_i)p^{(Z-\delta_i)/3}$

$= Z[(1 - \beta\omega)a^{3(1-\beta\omega)} + (1 + \beta\omega)a^{3(1+\beta\omega)} - a^6 - (1/3)\sum_{i=1}^{3}(1 - \delta_i\omega)a^{2(1-\delta_i\omega)}]$

$= Za^2[(1 - \beta\omega)a^{1-3\beta\omega} + (1 + \beta\omega)a^{1+3\beta\omega} - a^4 - (1/3)\sum_{i=1}^{3}(1 - \delta_i\omega)a^{-2\delta_i\omega}] = Za^2(S - a^4 - T)$.

By the arithmetic-geometric mean inequality,

$$S \geq 2((1 - \beta^2\omega^2)a^2)^{1/2} = 2a(1 - \beta^2\omega^2)^{1/2} \geq 2a(1 - (1/7)^2)^{1/2} > 1.97948a = 1.56.$$

Since $a^4 = c^{2/3} < .40$, we have $S - a^4 > 1.16$. Numerical experiments indicate that $T \leq 1$ but to give a proof of this seems difficult, so we provide a weaker upper bound for $T$. It is clear that if $Z \equiv 0 \mod 3$, then $\delta_i = 0$ for $i = 1, 2, 3$, so $T = 1$. If $Z \equiv 1 \mod 3$, then $\delta_1 = -2, \delta_2 = \delta_3 = 1$. Hence $3T = 2(1 - \omega)a^{-2\omega} + (1 + 2\omega)a^{4\omega}$.

Once can easily check that if $\delta \, \varepsilon \, \{-1, -2, 1, 2\}$, then the mapping $\omega \to (1 - \delta\omega)a^{-2\delta\omega}$ is a decreasing function of $\omega$ if and only if $\delta > 0$. It follows that $2(1 - \omega)a^{-2\omega} \leq 2$, and $(1 + 2\omega)a^{4\omega} \leq 1.15$ (max at $\omega = 1/7$) , so $T \leq 3.15/3 = 1.05$. If $Z \equiv 2 \mod 3$, then $\delta_1 = \delta_2 = -1$ and $\delta_3 = 2$. Therefore, $3T = 2(1 + \omega)a^{2\omega} + (1 - 2\omega)a^{-4\omega}$. We have $2(1 + \omega)a^{2\omega} \leq 2.14$ (max at $\omega = 1/7$), and $(1 - 2\omega)a^{-4\omega} \leq 1$, so that $T \leq 3.14/3 < 1.05$.

Combining these upper bounds for $T$ with the lower bound for $S - a^4$, we have $S - a^4 - T > 1.16 - 1.05 > 0$. Hence, $p_{2,3} < p = c^{1/Z} \leq p_{2,1}$.

Again, using Conjecture A, I can show:

**Theorem 30.** *Suppose that the school is profitable. Then*

*(a) as W decreases, the optimal number of classes does not decrease, and*

*(b) there is a proper subset $\bar{L}$ of the profitable operating region, such that if $(p, W)$ is in $\bar{L}$ and $p$ increases, then the optimal number of classes does not increase.*

**Proof.** Let $M_o, q_k, r_k, f_k$, and $Q_k$ be as in Definition 2. Put

$$L(k) = \{(p, W) | f_{k+1}(p) \leq W < f_k(p) \text{ and } W < V(Q_k, p)/k\}, \quad (2 \leq k \leq M_o),$$

$$L(1) = \{(p, W) | f_2(p) \leq W \text{ and } W < V(Q_1, p)\}, \text{ and}$$

$$L(Z) = R(Z).$$

Claim: $L(k) \cap L(j) = \emptyset$ if $k \neq j$. Without loss of generality we can assume $k > j$. We have to prove this in two cases: (1) $j > 1$ and (2) $j = 1$. Suppose case (1) holds. If $(p, W)\varepsilon L(k) \cap L(j)$, then

$$f_k(p) > W \geq f_{k+1}(p)$$

and

$$f_j(p) > W \geq f_{j+1}(p).$$

Hence, $f_k(p) > f_{j+1}(p)$, $f_k \neq f_{j+1}$, $k > j+1$, and $f_j \neq f_{j+1}$. By Proposition 22(a,b) and Proposition 29(a), $p < p_{j+1,k} \leq p_{j+1,j}$. But then $f_{j+1}(p) > f_j(p)$, a contradiction.

Second assume case (2) holds, so that $k > j = 1$. If $(p, W) \varepsilon L(k) \cap L(1)$, then

$$f_k(p) > W \geq f_{k+1}(p)$$

and

$$V(Q_1, p) > W \geq f_2(p).$$

Hence, $f_k(p) > f_2(p)$, $f_k \neq f_2$, $k > 2$, and $V(Q_1,.) \neq f_2$. By Proposition 22(a,b) and Proposition 29(b), $p < p_{2,k} \leq p_{2,1}$. But then $f_2(p) > V(Q_1, p)$, a contradiction.

Claim: $L(k) = R(k)$ whenever $1 \leq k \leq M_o$ or $k = Z$. If $k = Z$, this is obvious, so assume $k < Z$. Let $(p, W) \varepsilon R(k)$. By Theorem 18 the optimal class size vector is $Q_k$ and by Theorem 28, $k \leq M_o$. Hence,

$$V(Q_k, p) - kW > V(Q_{k-1}, p) - (k-1)W \quad (2 \leq k \leq M_o),$$

$$V(Q_k, p) - kW \geq V(Q_{k+1}, p) - (k+1)W \quad (1 \leq k \leq M_o),$$

and

$$V(Q_k, p) - kW > 0 \quad (1 \leq k \leq M_o).$$

Therefore

$$W < V(Q_k, p)/k \quad (1 \leq k \leq M_o),$$
$$f_k(p) = V(Q_k, p) - V(Q_{k-1}, p) > W \quad (2 \leq k \leq M_o),$$

$$W \geq V(Q_{k+1}, p) - V(Q_k, p) = f_{k+1}(p) \quad (1 \leq k \leq M_o),$$

so that $(p, W)$ is in $L(k)$. This shows that $R(k) \subseteq L(k)$. Since the $L(k)$ are disjoint sets and the union of the sets $R(k)$ is the entire region of profitabilty, it must be that $R(k) = L(k)$.

(a) If $p \leq 1/2$, then by Proposition 27 the optimal number of classes when the teacher cost is $W$ is $Z$. Hence, as $W$ decreases the number of classes does not decrease. Therefore, we can consider the case $p > 1/2$. Since the school is profitable, $(p, W) \varepsilon L(m)$, for some integer $m$. If $m = 1$, the result is clear, so assume $m > 1$. Then $f_{m+1}(p) \leq W < f_m(p)$. Let $N(k) = \{(p, W) | f_{k+1}(p) \leq W < f_k(p)\}$. Let $W' < W$. Then $W' < W \leq f_m(p)$. Suppose that $(p, W')$ is not in $N(k)$ for all $k \geq m$. The following is true if $k = m$:

$$W' < f_k(p). \tag{10}$$

If $W' \geq f_{k+1}(p)$, then $(p, W') \varepsilon N(k)$, a contradiction. Hence (10) holds with $k$ replaced by $k+1$ and we can continue the process. Thus, $W' < f_k(p)$, for all $k \geq m$. In particular, $0 < W' < f_{M_o+1}(p) = 2p(1-p)$. Since $p > 1/2$, $(p, W')\varepsilon L(Z)$, a contradiction. Hence, $(p, W')\varepsilon N(k)$ for some integer $k \geq m$. The school is profitable at $(p, W')$ since $V(Q_m, p) - mW' > V(Q_m, p) - mW > 0$. Hence, $(p, W')$ belongs to the intersection of $N(k)$ and the profitable region. But this intersection is $L(k)$. Since $k \geq m$ the result follows.

(b) Let $2 \leq k \leq M_o + 1$. Choose $s_k$ according to Proposition 22(d). Put $X_k = (s_{k+1}, f_{k+1}(s_{k+1}))$. The curve defined by $W = V(Q_k, p)/k$ forms part of the boundary of $L(k)$. If it meets the curve $W = f_{k+1}(p)$ to the left of $X_k$, then we define $L'(k)$ to be $L(k)$ with the small region consisting of the points of $L(k)$ that are to the left of $X_k$ and below the horizontal line $W = f_{k+1}(s_{k+1})$ removed. Otherwise, $L'(k)$ is just $L(k)$. Define $\bar{L} = \bigcup L'(k)$. The "profitability constraint" is defined as follows: If $W \geq f_2(s_2)$, it is the curve $W = V(Q_1, p)$. If $W \leq 2p(1-p)$, it is the straight line segment $W = p$. If $2p(1-p) < W < f_2(s_2)$, it is $\bigcup_{k=2}^{M_o}\{(p, V(Q_k, p)/k) | (p+\varepsilon, V(Q_k, p)/k) \text{ is in } L(k) \text{ for arbitrarily small positive } \varepsilon\}$. Finally, add to the above points all limit points. We have to show that the piecepartsfit together so that we have a function of $p$. To do this we have to show that the curves $W = V(Q_k, p)/k$ and $W = V(Q_{k-1}, p)/(k-1)$ meet $W = f_k(p)$ at the same point. Suppose $f_k(p_o) = V(Q_k, p_o)/k$. Then $V(Q_k, p_o) - V(Q_{k-1}, p_o) = V(Q_k, p_o)/k$, so that $((k-1)/k)V(Q_k, p_o) = V(Q_{k-1}, p_o)$, or $V(Q_{k-1}, p_o)/(k-1) = V(Q_k, p_o)/k = f_k(p_o)$. Since each $V(Q_k, p)$ is an increasing function of $p$, it is clear that the profitability constraint is an increasing function of $p$. We denote this constraint function by $C(p)$.

The lower boundary of $L'(k)$ is defined by a non-increasing function of $p$, which we denote by $\tilde{f}_{k+1}(p)$. Note that $\tilde{f}_{k+1}(p) \geq f_{k+1}(p)$. Suppose $(p, W)\varepsilon\bar{L}$, say $(p, W)\varepsilon L'(m)$. If $m = Z$ there is nothing to show, so assume $m < Z$. Then by Theorem 28 $m \leq M_o$. There are two cases: (1) $m = 1$, and (2) $2 \leq m \leq M_o$. Let $p' > p$. First, consider case (1). Here $(p, W)\varepsilon L'(1)$, so $\tilde{f}_2(p) \leq W < V(Q_1, p)$. Hence $W \geq \tilde{f}_2(p) \geq \tilde{f}_2(p')$. And $W < V(Q_1, p) < V(Q_1, p')$, so $(p', W)\varepsilon L'(1)$ and the optimal number of classes has not increased. Second, consider case (2). Let $\tilde{N}(k) = \{(p, W) | \tilde{f}_{k+1}(p) \leq W < \tilde{f}_k(p)\}$, and $\tilde{N}(1) = \{(p, W) | \tilde{f}_2(p) \leq W < C(p)\}$. Since $(p, W)\varepsilon L'(m)$, we have $W < C(p)$ and $\tilde{f}_{m+1}(p') \leq \tilde{f}_{m+1}(p) \leq W < \tilde{f}_m(p)$. Suppose that $(p', W)$ is not in $\tilde{N}(k)$ for all $k \leq m$. Then the following is true for $k = m$:

$$W \geq \tilde{f}_k(p'). \tag{11}$$

Since $(p', W)$ is not in $\tilde{N}(k)$, it must be that $W \geq \tilde{f}_{k-1}(p')$, that is, (11) holds with $k$ replaced by $k - 1$.

Continuing this process we get $W \geq \tilde{f}_2(p')$ and $W \geq C(p')$. Hence, $W \geq C(p') > C(p) > W$, a contradiction. Therefore, $(p', W) \varepsilon \tilde{N}(k)$ for some $k \leq m$. If the profitability constraint function $C$ intersects $\tilde{f}_k$ to the left of $X_{k-1}$, then denote by $T$ the small region that is excluded when defining $L'(k-1)$. Suppose that $(p', W) \varepsilon T$. Then $f_k(p') \leq W < \tilde{f}_k(s_k)$ and $p' < s_k$. Hence, $f_k(p) < f_k(p') \leq W < f_k(s_k)$. and $(p, W) \varepsilon T$, a contradiction. Therefore, $(p', W)$ is not in $T$. It follows that $(p', W) \varepsilon L'(k)$, so that the number of classes has not increased.

For $Z = 5$ there is a small region $\Gamma$ within the profitable operating region such that if $(p, W)$ is in $\Gamma$, and if $p$ increases by a sufficiently large amount, then $m$ also increases. Take $W = .673$. One can check that for $p = .60$ the school is profitable and the optimal class size vector is (2,3). But for $p = .62$ the optimal class size vector is (1,2,2), so that $m$ has increased. This shows that only a proper subset of the profitable operating region will have the property of decreasing $m$ with increasing $p$.

It can be shown that for $Z \geq 3$ and odd, there is always a small region $T \subseteq L(1)$ such that if $(p, W) \varepsilon T$ then there is a $p' > p$ with $(p', W)$ in $L(2)$, so that the optimal number of classes increases with $p$. The reason for this is that the profitabilty constraint meets $f_2$ to the left of the maximum point $X_1$.

## 5. More than one type of student

In this section I extend Lazear's basic model to cover the case where there are $s > 1$ types of students.

I assume that there are $s$ different types of students, say $1, 2, \ldots, s$ and they have probabilities $p_1, \ldots, p_s$ of not disrupting the class. Let there be $a_i$ students of type $i$, so that $\sum_{i=1}^{s} a_i = Z$. The school seeks to maximize

$$\sum_{j=1}^{m} (\sum_{i=1}^{s} n_{i,j}) \prod_{i=1}^{s} p_i^{n_{i,j}} - Wm$$

subject to

$$\sum_{j=1}^{m} n_{i,j} = a_i \ (1 \leq i \leq s), \ 0 \leq n_{i,j} \ \text{and} \ n_{i,j} \ \text{an integer.}$$

In general there will be some mixed classes having students of different types and some segregated classes having students all of a single type. While Lazear argued that segregation was optimal, this is not correct.

In general there will be some mixing. For example if $Z = 6, a_1 = a_2 = 3, p_1 = .8, p_2 = .5$ and $W = .51$, then the maximum profit of 1.05 is attained at a class size vector of (2,2,2), with the type 1's distributed as (0,2,1) and the type 2's as (2,0,1). Thus, the third class has a mix of type 1 and type 2 students.

It is of interest to understand therefore the structure of an optimal solution, in particular, how many mixed classes there can be and if there is any pattern to the mixing. Within the set of segregated classes the distribution of students and number of classes is determined by Theorem 18: the class sizes are nearly equal. So we focus on the mixed classes.

**Theorem 31**. *If p and q are distinct positive reals and u,v are positive integers, then* $(p^{u-1}q^{v+1} + p^{u+1}q^{v-1})/2 > p^u q^v$.

**Proof**. This follows immediately from the arithmetic-geometric mean inequality.

We assume that $\{n_{i,j}\}$ is an optimal solution and that there are $m$ classes in the solution. Let $A = (n_{i,j})$ be the matrix that describes the distribution of students within classes. $G(A)$ is the corresponding bipartite graph of $A$. That is, $G(A) = (R, S, E)$ has "row" nodes $R = \{1, 2, \ldots, s\}$ and "column" nodes $S = \{1, \ldots, m\}$. Its edges are $E = \{(i, j) | n_{i,j} > 0\}$. $\tilde{A}$ will denote the submatrix of $A$ whose rows are $1, \ldots, s$ and whose columns correspond to the mixed classes.

**Theorem 32**. *$G(A)$ is a forest and the number of mixed classes is at most s-1. The number of mixed classes is exactly s-1 if and only if $G(\tilde{A})$ is a path with $2(s-1)$ edges.*

**Proof.** Suppose that $G(A)$ has a cycle of length $2t$. Then after permutation of the rows and columns of $A$ we can assume that $n_{1,1}, n_{2,1}, n_{2,2}, \ldots, n_{t,t}, n_{1,t}$ are all positive. Moreover, because each of these entries is positive, each is also less than the bound for its row. Let $\Delta = (\delta_{i,j})$ be the $s \times m$ matrix with $\delta_{i,i} = 1$, if $1 \leq i \leq t$, $\delta_{i+1,i} = -1$ if $1 \leq i \leq t-1$ and and $\delta_{1,t} = -1$. All other entries are 0. Then $A + \Delta$ and $A - \Delta$ are both feasible solutions for the maximization problem. Let

$$f(A) = \sum_{j=1}^{m}(\sum_{i=1}^{s} n_{i,j})\prod_{i=1}^{s} p_i^{n_{i,j}}.$$

Put $N_j = \sum_{i=1}^{s} n_{i,j}$ and $C_j = \prod_{\substack{i=1 \\ i \neq j, j+1}}^{s} p_i^{n_{i,j}}$, where $j+1$ is always taken mod $t$. Then for $\varepsilon = 1$ or $\varepsilon = -1$, we have

$f(A + \varepsilon\Delta) = \sum_{j=1}^{m} N_j C_j p_j^{n_{j,j}+\varepsilon\delta_{j,j}} p_{j+1}^{n_{j+1,j}+\varepsilon\delta_{j+1,j}}$, so that by Theorem 31

$$(f(A+\Delta) + f(A-\Delta))/2 = \sum_{j=1}^{m} N_j C_j (p_j^{n_{j,j}+\delta_{j,j}} p_{j+1}^{n_{j+1,j}+\delta_{j+1,j}} + p_j^{n_{j,j}-\delta_{j,j}} p_{j+1}^{n_{j+1,j}-\delta_{j+1,j}})/2$$
$$> \sum_{j=1}^{m} N_j C_j p_j^{n_{j,j}} p_{j+1}^{n_{j+1,j}}$$
$$= f(A).$$

Hence, either $f(A+\Delta) > f(A)$ or $f(A-\Delta) > f(A)$ and $A$ cannot be a maximum for $f$, a contradiction. Therefore, $G(A)$ has no cycles and is a forest. Since $G(\tilde{A})$ is a subgraph of $G(A)$ it is a forest. Let $\tilde{m}$ be the number of mixed classes. Let $t$ be the number of trees in $G(\tilde{A})$. Since $G(\tilde{A})$ has $s + \tilde{m}$ nodes it follows that $G(\tilde{A})$ has $s + \tilde{m} - t$ edges. But each column node corresponding to a mixed class has degree at least 2, so $G(\tilde{A})$ has at least $2\tilde{m}$ edges. Hence $s + \tilde{m} - t \geq 2\tilde{m}$, and $\tilde{m} \leq s - t \leq s - 1$. If equality holds then $t = 1$ and each column node for the mixed classes has degree 2. This implies that $G(\tilde{A})$ is a path with $2(s-1)$ edges.

Some examples based on numerical investigations of the case of $s = 4$ reveal some interesting differences between the $s = 1$ and $s > 1$ cases.

We represent the problem and its solution by an $s$-rowed matrix, where the first column has the values of $p_1, \ldots, p_s$ and the second column has the values of $a_1, \ldots, a_s$. The remaining columns represent the class sizes with the mixed classes presented first. In all the examples $W$ is .9.

Class sizes are for the segregated class sizes are not necessarily nearly equal when $s > 1$.

| .5  | 13 | 0 | 13 |   |   |   |
|-----|----|---|----|---|---|---|
| .75 | 5  | 2 |    | 3 |   |   |
| .8  | 4  | 1 |    |   | 3 |   |
| .85 | 3  | 0 |    |   |   | 3 |

Class sizes for the mixed classes may not be nearly equal and $G(A)$ may be 2 disjoint trees.

| 1 | 3 | 0 | 3 |   |   |
|---|---|---|---|---|---|
| .81 | 4 | 1 | 0 | 3 |   |
| .82 | 4 | 0 | 1 |   | 3 |
| .707 | 3 | 1 | 0 |   | 2 |

The first two columns are for the mixed classes which have sizes 2 and 4.

There may be $s - 1$ mixed classes.

| .8 | 3 | 0 | 1 | 2 |   |   |
|---|---|---|---|---|---|---|
| .81 | 4 | 2 | 2 | 0 |   |   |
| .82 | 4 | 1 | 0 | 0 | 3 |   |
| .707 | 3 | 0 | 0 | 1 | 0 | 2 |

Proposition 26 has a generalization for the case of $s > 1$ types of students.

**Theorem 33.** *Suppose that* $W \geq \max\{p_i + p_j - 2p_i p_j\}$. *Then an optimal solution has at most one singleton class.*

**Proof.** If there are two singleton classes, say with one of type $i$ and the other of type $j$ students, the profit from these two classes is $p_i + p_j - 2W$. This profit should exceed the profit from combining the classes, which is $2p_i p_j - W$ and that implies $p_i + p_j - 2p_i p_j > W$, a contradiction.

Using Theorem 33 I can now prove that there are gaps in the possible number of classes for an optimum solution.

**Theorem 34**. *Suppose that* $W \geq \max\{p_i + p_j - 2p_i p_j\}$. *For a profitable school of Z students the optimal number of classes m satisfies m = Z , or m $\leq$ Z/2, Z even or m $\leq$ (Z + 1)/2, Z odd.*

**Proof.** By Theorem 33 there is a most one singleton class so the other classes have size at least 2. As in the proof of Theorem 28, it follows that $m \leq Z/2$ or $m \leq (Z + 1)/2$.